\documentclass{ifacconf}

\usepackage{graphicx}      % include this line if your document contains figures
\usepackage{natbib}        % required for bibliography

\usepackage{amsmath,amssymb,amsfonts}
\usepackage{textcomp}
\usepackage{mathtools}

\usepackage{amsmath} % assumes amsmath package installed
\usepackage{amssymb}  % assumes amsmath package installed

\usepackage{accents}
\newcommand{\ubar}[1]{\underaccent{\bar}{#1}}

\newtheorem{definition}{Definition}[section]
\newtheorem{assumption}[definition]{Assumption}

\newtheorem{remark}[definition]{Remark}

\newtheorem{problem}[definition]{Problem}

\DeclareMathOperator*{\minimize}{\ \ minimize\ \ }
\DeclareMathOperator{\st}{\ \ subject\ to\ }

\DeclareMathOperator{\CoP}{CoP}
\DeclareMathOperator{\cpu}{cpu}

\allowdisplaybreaks

\DeclareSymbolFont{bbold}{U}{bbold}{m}{n}
\DeclareSymbolFontAlphabet{\mathbbold}{bbold}

\usepackage{calrsfs}
\DeclareMathAlphabet{\pazocal}{OMS}{zplm}{m}{n}
\renewcommand{\mathcal}[1]{\pazocal{#1}}

\usepackage{siunitx}

\usepackage[usenames]{xcolor}

\usepackage{mathtools}

\begin{document}
\begin{frontmatter}

\title{Model Predictive Control for Energy-Efficient Operations of Data Centers with Cold Aisle Containments\thanksref{footnoteinfo}} 

\thanks[footnoteinfo]{This work is funded in part by JSPS KAKENHI Grant Number 18K13777, the National Natural Science Foundation of China (NSFC) under Grant No.~61502255, and Inner Mongolia Science and Technology Plan.}

\author[naist]{Masaki Ogura} 
\author[Wan]{Jianxiong Wan} 
\author[naist]{Shoji Kasahara}

\address[naist]{Nara Institute of Science and Technology, 
Ikoma, Nara 630-0192 Japan (e-mail: 
oguram@is.naist.jp, kasahara@is.naist.jp)}
\address[Wan]{Inner Mongolia University of Technology,
Hohhot, 010080, Inner Mongolia, China (e-mail: jxwan@imut.edu.cn)}
%\address[Third]{Electrical Engineering Department, 
%   Seoul National University, Seoul, Korea, (e-mail: author@snu.ac.kr)}

\begin{abstract}
In this paper, we study a problem of controlling cooling facilities and computational equipments for energy-efficient operations of data centers. Although a plethora of approaches have been proposed in previous literatures, there is a lack of rigorous methodologies for effectively addressing the problem. To bridge this gap, we propose MPC frameworks for jointly and optimally tuning the Computer Room Air Conditioner (CRAC) supplying air temperature as well as the number of active servers for minimizing the overall energy consumption. We specifically find that, when a standard model of data centers with contained cold aisles and cooling facilities are used, the optimization problems arising from the MPC framework can be transformed to  convex optimization problems that can be solved efficiently. We present several numerical simulations to illustrate the effectiveness of our theoretical results.
\end{abstract}


\end{frontmatter}
%===============================================================================

%%%%%%%%%%%%%%%%%%%%%%%%%%%%%%%%%%%%%%%%%%%%%%%%%%%%%%%%%%%%%%%%%%%%%%%%%%%%%%%%

\section{Introduction} 

Data centers are major energy consumers in the current Information Technology (IT) industry. It is predicted that data centers in the U.S.~alone will consume $140$ billion~kWH electricity by~$2020$~\citep{Delforge2014}. The huge energy consumption not only poses an enormous monetary burden to data center operators, but also leads to serious concerns for environmental issues~\citep{Gao2012b}. Therefore, it is imperative for researchers to find efficient ways to cut down the data center energy consumption.

Dissecting the energy consumption of a data center shows that the IT and cooling subsystems play dominant roles~\citep{Dayarathna2016}. Unfortunately, majority of the works from the Computer Science community concentrated only on the control of IT facilities. \cite{Lin2016a} used an integer linear programming to model the energy minimization problem where the control variables were the state of servers (on/off) and workload dispatching policy. An enhanced genetic algorithm and a heuristic greedy sequence approach were proposed to solve the integer linear programs. VMAP+ proposed by~\cite{Lee2017b} reduces the energy consumption via allocation and migration of virtual machines. 
%SpAWM \citep{Chen2011a} and ThermoRing \citep{Zhao2016b} are workload dispatchers which try to raise the server inlet temperature by optimizing workload/heat distribution so that the cooling energy is reduced. 
TIGER~\citep{Chavan2016a} also uses the same idea to cut the cooling cost in storage clusters but the control knob is the file placement strategy. Although the above IT side control methods are thermal-aware to some extent, they can only affect the heat generation pattern of the data center and have limited impact on the cooling subsystem.

Recently, joint optimization techniques which integrate control knobs from both IT and cooling subsystems are proposed to balance the energy consumption in these components. \cite{Parolini2012} introduced a control framework which includes both thermal management and workload scheduling. Under the simplified assumption that the temperature of rack exhaust air is constant, \cite{Li2012e} derived a closed form solution to the joint optimization problem. \cite{Fang2016} used the framework of Model Predictive Control (MPC) to coordinate the actions of Computer Room Air Conditioners (CRACs) and servers. %They further proposed a two-time scale control framework where workload dispatching and CPU frequency were adjusted on a coarse time scale and CRAC supplying air temperature was manipulated in a fine time scale \citep{Fang2017}.

In this paper, we present MPC formulations for achieving overall energy minimization in data centers. We specifically study the data centers with Cold Aisle Containments (CACs) that are adopted in nearly all modern data centers. Furthermore, we show that the finite-horizon control problems that have to be iteratively solved under MPC formalisms can be transformed to convex optimization problems, which allows us to efficiently find the optimal control inputs by off-the-shelf solvers. These features are in contrast to the previous works in \cite{Parolini2012,Fang2016}, where data centers with open aisle systems involving hot air recirculations are considered while the optimality of the control inputs is not guaranteed. 

This paper is organized as follows. In Section~\ref{sec:problemFormulation}, we introduce the model of data centers and formulate the energy minimization problem. In Section~\ref{sec:deterministic}, the energy minimization problem is cast into an MPC problem under the assumption that user request rates are predictable. In Section~\ref{sec:stochastic}, we present a scenario-based MPC formulation in which we do not require predictions but utilize a history of past data. We illustrate the effectiveness of the proposed methods via numerical simulations in Section~\ref{sec:simulation}.

\section{Problem Formulation}\label{sec:problemFormulation}

We formulate the problem studied in this paper. We first present the model of data centers with CACs in Section~\ref{ref:model}. In  Section~\ref{sec:problem}, we formulate the overall energy minimization problem.

\subsection{Data Center Description}\label{ref:model}

\begin{figure}[tb]
\centering
\includegraphics[width =.95\linewidth]{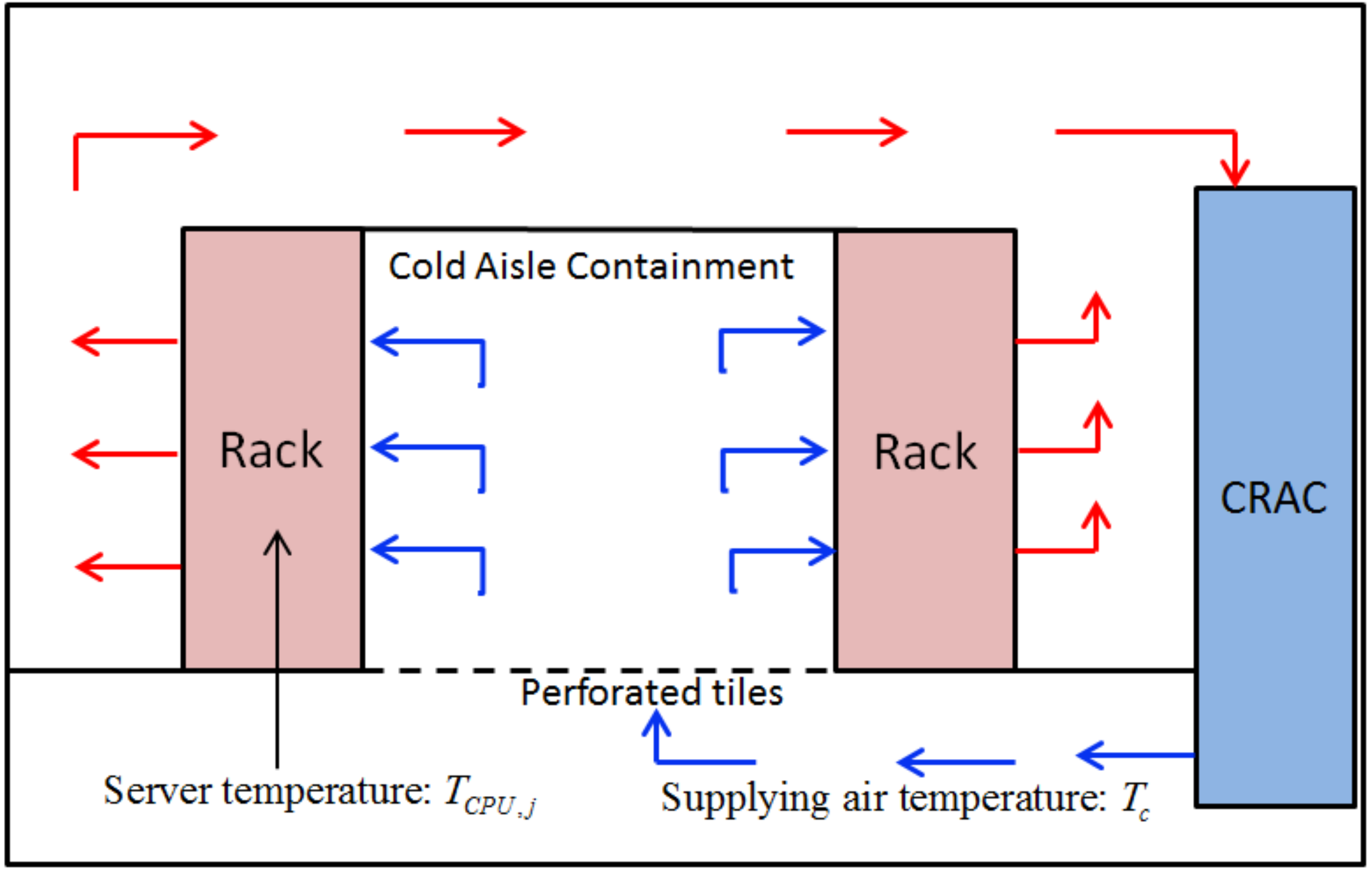}
\caption{A typical machine room layout with cold aisle containment inside data centers.}
\label{fig:DCMR}
\end{figure}

Let us first briefly overview the structure of the data centers studied in this paper. A typical machine room layout for modern data centers is shown in Fig. \ref{fig:DCMR}. Racks are placed on the raised floor with inlets faced towards each other to form a cold aisle. Almost all power drawn by servers is dissipated as heat, which jeopardizes the reliability of servers. A CRAC is installed near the side wall of the machine room and constantly blows cold air (blue arrows) into the underfloor plenum, leading to a positive pressure differential between the plenum and cold aisle which forces the cold air to be ejected from the perforated tiles located in front of racks. With the help of built-in fans, servers suck in the cold air, dump the waste heat into it, and finally discharge it into the hot aisle. The hot airflow (red arrows) will return to the CRAC though the ceiling.

In traditional open aisle data centers, the hot air recirculation, which is a phenomenon that the hot air reenters the cold aisle and mixes with the supplying cold air, can be easily observed in the rack top and aisle ends due to the nonuniform pressure distribution inside the cold aisle. Since the hot air recirculation compromises the cooling efficiency, researchers developed containment systems for cold or hot aisles. The containment systems usually consist of concrete barrier plates installed at rack top and aisle ends which minimize the effect of air mixing. It was shown that CACs provide ``close to perfect cold air delivery''~\citep{Arghode2013}.

We then describe the mathematical model of the data centers \citep{Fu2017,Parolini2012}. We assume that there are a total number of~$J$ users requesting service from a data center. For each nonnegative integer~$t$, let $L_j(t)$ denote the request rate of user~$j$ at time~$t$. The request from user~$j$ is assumed to be processed in the $j$th cluster of servers located in the data center. We assume that the number of the servers in the $j$th cluster, denoted by~$m_j(t)$, can be dynamically tuned by the operator of the data center. We denote  by~$T_{\cpu, j}(t)$ the temperature of the CPUs in the $j$th cluster at time~$t$.

We adopt a standard model of the CPU temperatures, which is described below. We model the energy consumption of a single server in the $j$th cluster by
\begin{equation}\label{eq:def:p_j}
p_j(t) = a_1 \frac{L_j(t)}{m_j(t)} + a_2, 
\end{equation}
where $a_1$ is the marginal energy consumption for increased CPU utilization, and $a_2$ is the energy consumption of remaining components apart from CPU \citep{Dayarathna2016}. Therefore, these two parameters are all positive constants. Then, the dynamics of the CPU temperature is described by the linear difference equation
\begin{equation}\label{eq:temperatureDynamics}
T_{\cpu, j}(t+1) = \alpha_j T_c(t) +  \varsigma_j p_j(t) + \beta_j T_{\cpu, j}(t), 
\end{equation}
where $T_c(t)$ is the CRAC supplying air temperature and $\alpha_j, \beta_j, \varsigma_j$ are heat exchange rates~\citep{Parolini2012}.

\subsection{Problem Formulation}\label{sec:problem}

The energy consumption of a server cluster handling the requests from user~$j$ at time~$t$ is modeled as $P_j(t) = p_j(t)m_j(t)$. Then, the total energy consumption of the data center at time~$t$ can be expressed by $P(t) = \sum_{j=1}^J P_j(t)$. On the other hand, the energy consumption at the CRAC is modeled by 
\begin{equation}\label{eq:def:C}
C(t) = P(t)/\CoP(T_c(t)), 
\end{equation}
where $\CoP(T) = 0.0068T^2 + 0.0008T + 0.458$ is a function describing the cooling efficiency of CRAC at the temperature~$T$ and is widely adopted in the literature~\citep{Moore2005}. Therefore, the total energy consumption in the data center at time~$t$ equals $E(t) = P(t) + C(t)$.

We finally state the constraints that have to be satisfied during the operation of the data center. First, the temperature of the CRAC can be controlled only within the following interval
\begin{equation}\label{eq:const:T_c}
\ubar T_c \leq T_c(t)\leq \bar T_c, 
\end{equation}
where $\ubar T_c$ and $\bar T_c$ are given positive constants. The second requirement is on the response times of users. Let $D_j(t)$ denote the response time of user $j$. We require that
\begin{equation}\label{eq:delayConst}
D_j(t) \leq \bar D_j, 
\end{equation}
where $\bar D_j$ is an upper bound of tolerable response time specified by user $j$. We assume the system is a classical $M/M/1$ queueing system~\citep{Harchol-Balter2013}. Under this assumption, the response time admits the representation~$D_j(t) = {1}/{(m_j(t)\mu_j - L_j(t))}$, where $\mu_j$ is the average service rate for user~$j$. 

We finally require that the CPU temperatures are kept below a certain safety level. Specifically, we require that
\begin{equation}\label{eq:cpuTemperatureConst}
T_{\cpu, j}(t) \leq \bar T_{\cpu, j},
\end{equation} 
where $\bar T_{\cpu, j} > 0$ is a constant given by hardware manufacturers to ensure the server reliability. 

We can now state our energy minimization problem studied in this paper: 

\begin{problem}\label{prb:dyn}
Find a sequence of the CRAC temperatures~$T_c = \{T_c(t)\}_{t=t_0}^{t_f}$ and the amount of active servers~$m_j = \{m_j(t)\}_{t=t_0}^{t_f}$ ($j=1$, \dots, $J$) solving the following optimization problem:
\begin{subequations}\label{eq:opt:original}
\begin{align}
\minimize_{T_c,\,m_1,\,\dotsc,\,m_J}\ \ \ & \sum_{t=t_0}^{t_f} E(t)\label{eq:objective}
\\
\st\ \ \ & 
%\mbox{\eqref{eq:const:T_c}, \eqref{eq:delayConst},  \eqref{eq:cpuTemperatureConst}}.  
\mbox{\eqref{eq:const:T_c}--\eqref{eq:cpuTemperatureConst}}.  
\end{align}
\end{subequations}
\end{problem}

The difficulty of Problem~\ref{prb:dyn} mainly stems from the nonlinearity of the objective function~\eqref{eq:objective} as well as the constraint~\eqref{eq:cpuTemperatureConst}. The energy consumption~\eqref{eq:def:C} at CRAC involves the inverse of a quadratic function of the cooling temperature. Furthermore, by equations~\eqref{eq:def:p_j} and~\eqref{eq:temperatureDynamics}, the CPU temperature involves the inverse of $m_j$. For these reasons, it is not a trivial task to directly solve the optimization problem~\eqref{eq:opt:original}. 
%Another difficulty lies in the availability of the number of user requests. To solve Problem~\ref{prb:dyn}, we need, at time $t=t_0$, the values of $L_j(t)$ in future, i.e., for all $t_0\leq t\leq t_f$. However, such predictions are hardly precise in practice. 

\section{Model Predictive Control}\label{sec:deterministic}

In this section, we present an MPC formulation for solving problem \ref{prb:dyn}. Assuming that a prediction on the user request rate is available at each time instant, we formulate the energy minimization problem as an iteration of optimal control problems within smaller time-windows. We show that, under a mild and reasonable assumption on the CRAC supplying air temperature, the problem can be transformed to a convex optimization problem.

\subsection{Convexity}

In this section, we place the following assumption: 

\begin{assumption}\label{asm:deterministic}
There exists a positive integer~$t_h$ such that, at each time~$\tau \geq t_0$, we can predict the number of requests~$\{L_j(t)\}_{t=\tau}^{\tau + t_h}$ for all user~$j$.
\end{assumption}

This assumption allows us to formulate the following finite-horizon optimal control problem for $\tau = t_0, \dotsc, t_f$: 
\begin{subequations}\label{eq:mpc:deterministic}
\begin{align}
\minimize_{T_c,\,m_1,\,\dotsc,\,m_J,\,\gamma}\ \ \ & \gamma
\\
\st\ \ \ \,& \mbox{\eqref{eq:const:T_c}--\eqref{eq:cpuTemperatureConst}},
\\
& \sum_{t=\tau}^{\tau + t_h}  E(t) \leq \gamma,  \label{eq:performance_gamma}
\end{align}    
\end{subequations}
where we have introduced the slack variable~$\gamma$. 

In order to state the main result of this section, we introduce the following notations. Let $\tau \geq t_0$ be arbitrary. For $t \in \{\tau, \tau+1, \dotsc, \tau + t_h\}$, let $T_{\cpu, j} (t; u_\tau)$ denote the temperature~$T_{\cpu, j} (t)$ when we apply the control inputs
\begin{equation}\label{eq:def:realInputSeq}
\begin{multlined}[.9\linewidth]
u_\tau = \{T_c(t), m_1(t), \dotsc, m_J(t)\colon  \\ 1\leq j\leq J,\,\tau\leq t\leq \tau+t_h\} \in (\mathbb{R}_+^{t_h+1})^{J+1}, 
\end{multlined}
\end{equation}
where $\mathbb{R}_+$ denotes the set of positive numbers. The notation~$E(t; u_\tau)$ is understood in the same manner. The following theorem shows that we can transform the MPC problem~\eqref{eq:mpc:deterministic} to a convex optimization problem and is the first main result of this paper.

\begin{thm}\label{thm:deterministic}
Assume that
\begin{equation}\label{eq:asm:T>11}
\ubar T_c \geq 11. 
\end{equation}
Then, the solution of the optimization problem~\eqref{eq:mpc:deterministic} is given by
\begin{equation}\label{eq:optimalInputRepresentation}
T_c(t) = e^{x^\star(t)},\ m_j(t) = e^{y_j^\star(t)}, 
\end{equation}
where $x^\star = \{x^\star(t)\}_{t=\tau}^{\tau+t_h}$ and $y_j^\star = \{y_j^\star(t)\}_{t=\tau}^{\tau+t_h}$ ($j=1, \dotsc, J$) solve the convex optimization problem
\begin{subequations}\label{eq:convexProg:deterministic}
\begin{align}
\minimize_{x,\,y_1,\,\dotsc,\,y_J,\,\Gamma}\ \ \ & \Gamma 
\\
\st\ \ \,& u_\tau = 
\begin{multlined}[t]
\{e^{x(t)}, e^{y_1(t)}, \dotsc, e^{y_J(t)} \colon \\  1\leq j\leq J,\,\tau\leq t\leq \tau+t_h\},\end{multlined}\label{eq:exponentialInputSeq}
\\
&
\log \ubar{T}_c \leq x(t)\leq \log \bar T_c,  \label{eq:x_bound}
\\
& y_j(t) \geq \log(\bar D_j^{-1} + L_j(t)) - \log \mu_j, 
\label{eq:exp:delayConst}
\\
&
\log T_{\cpu, j}(t; u_\tau) \leq \log \bar T_{\cpu, j}, 
\label{eq:exp:cpuTempConst}
\\
& \log \sum_{t=\tau}^{\tau + t_h}  E(t;  u_\tau) \leq \Gamma. 
\label{eq:lessthanGamma}
\end{align}
\end{subequations}    

\end{thm}

\begin{remark}
%The above analysis indicates that problem \eqref{eq:convexProg:deterministic} is a convex optimization problem under constraint \eqref{eq:asm:T>11}.
As can be seen in ASHRAE's Thermal Guidelines \citep{ASHRAE2011}, data center operators usually choose a supplying air temperature between $18$ and $27^\circ$C. Therefore, in most of the standard cases, the optimization problem~\eqref{eq:convexProg:deterministic} has a unique minimum which can be obtained by traditional convex optimization techniques.
\end{remark}

In the proof of Theorem~\ref{thm:deterministic}, the following lemma plays an important role. 

\begin{lem}\label{lem:convexity:1+1/cop}
Define the function $f \colon (0, \infty) \to (0, \infty)$ by 
\begin{equation}\label{eq:def:fcop}
f(T) = 1 + \CoP(T)^{-1}. 
\end{equation}
Then, the function 
\begin{equation}\label{eq:def:F}
F\colon \mathbb{R} \to \mathbb{R} \colon x \mapsto \log(f(e^x))
\end{equation}
is convex on the interval $[\log(11), \infty)$. 
\end{lem}

\begin{pf}
A straightforward calculation shows that 
there exists a polynomial~$p$ having positive coefficients such that
%\begin{equation*}
$F''(\log T) = T{p(T)^{-2}}q(T)$, 
%\end{equation*}
where $q(T) = q_5 T^5 + q_4 T^4 + q^3T^3 + q^2 T^2 + q_1T + q_0$ is the polynomial with the coefficients $q_5 = \num[group-separator={,}]{9826}$, $q_4 =\num[group-separator={,}]{2023}$, $q_3 =136$, $q_2 =-\num[group-separator={,}]{81426}$, $q_1 =-\num[group-separator={,}]{141899850}$, and $q_0 =-\num[group-separator={,}]{4173525}$. The polynomial~$q$ has the real zeros at $T = -10.99$, $-0.029$, and $10.94$. Furthermore, $q(11)  = \num[group-separator={,}]{37362464}$ is positive. Therefore, by the continuity of polynomials, the derivative~$F''(\log T)$ is positive for all $T \in [11, \infty)$.
\end{pf}

We also review the notion of posynomials~\citep{Boyd2007} for the proof of Theorem~\ref{thm:deterministic}. Let $f \colon \mathbb{R}^n_{+}\to \mathbb{R}_{+}$ be a function. We say that $f$ is a \emph{monomial} if there exist $c>0$ and real numbers~$a_1$, \dots,~$a_n$ such that $f(v) = c v_{\mathstrut 1}^{a_{1}} \dotsm v_{\mathstrut n}^{a_n}$. We say that $f$ is a \emph{posynomial} if $f$ is the sum of finitely many monomials. The following lemma shows the log-log convexity of posynomials and is used in the proof of Theorem~\ref{thm:deterministic}. 

\begin{lem}[\cite{Boyd2007}]\label{lem:llconvexity}
Let $f \colon \mathbb{R}_{+}^n \to \mathbb{R}_{+}$ be a posynomial. Define the function $F\colon \mathbb{R}^n \to \mathbb{R}$ by $F(w) = \log f(\exp[w])$, where $\exp[\cdot]$ denotes the entry-wise exponentiation. Then, $F$ is convex.
\end{lem}

Let us prove Theorem~\ref{thm:deterministic}. 

\begin{pf}
Under the transformations $x(t) = \log T_c(t)$, $y_j(t) = \log m_j(t)$, and $\Gamma = \log \gamma$, the constraints~\eqref{eq:const:T_c}, \eqref{eq:delayConst}, \eqref{eq:cpuTemperatureConst}, and \eqref{eq:performance_gamma} are equivalent to the constraints~\eqref{eq:x_bound}, \eqref{eq:exp:delayConst}, \eqref{eq:exp:cpuTempConst}, and \eqref{eq:lessthanGamma}, respectively. Therefore, the solution of the optimization problem~\eqref{eq:mpc:deterministic} is given by \eqref{eq:optimalInputRepresentation} with $x^\star$, $y_1^\star$, \dots, $y_J^\star$ being the solutions of the optimization problem~\eqref{eq:convexProg:deterministic}.

Let us show the convexity of the optimization problem~\eqref{eq:convexProg:deterministic} under the assumption~\eqref{eq:asm:T>11}. It is sufficient to show the convexity of the mappings $\psi \colon %(\mathbb{R}^{t_h+1})^{J+1} \to \mathbb{R} \colon
(x, y_1, \dotsc, y_J) \mapsto \log T_{\cpu, j}(t; u_\tau)$ and $\phi
%\colon (\mathbb{R}^{t_h+1})^{J+1} \to \mathbb{R} 
\colon (x, y_1, \dotsc, y_J) \mapsto \log \sum_{t=\tau}^{\tau + t_h} E(t; u_\tau)$. Equation \eqref{eq:temperatureDynamics} implies that $T_{\cpu, j}(t; u_\tau)$ is a posynomial in the variables in~\eqref{eq:def:realInputSeq} because 
\begin{equation*}
\begin{multlined}[.9\linewidth]
T_{\cpu, j}(t; u_\tau) = \beta_j^{t-\tau} T_{\cpu, j}(\tau) \\+ \sum_{s = 0}^{t-\tau} \beta_j^{t-\tau-s} \biggl(\alpha_j T_c(s) + a_1\varsigma_j \frac{L_j(s)}{m_j(s)} + a_2\varsigma_j\biggr). 
\end{multlined}
\end{equation*}
Therefore, Lemma~\ref{lem:llconvexity} implies the convexity of~$\psi$. Also, Equation~\eqref{eq:def:C} shows $E(t;u_\tau) = P(t;u_\tau) f(e^{x(t)})$ for the function~$f$ defined in \eqref{eq:def:fcop}.
Therefore, to show the convexity of~$\phi$, it is sufficient to show that the function~$F$ defined in~\eqref{eq:def:F} and the mapping $\theta\colon 
%(\mathbb{R}^{t_h+1})^{J+1} \to \mathbb{R} \colon
(x, y_1, \dotsc, y_J) \mapsto \log P(t; u_\tau)$ are both convex. The convexity of the mapping~\eqref{eq:def:F}  follows from Lemma~\ref{lem:convexity:1+1/cop} and the assumption~\eqref{eq:asm:T>11}. We can prove the convexity of $\theta$ in a similar way as $\psi$. We omit the details of the proof due to limitations of space. 
\end{pf}

\subsection{Suppression of Fluctuations}

As we observe in Section~\ref{sec:simulation}, the input sequence obtained by iteratively solving the optimization problem~\eqref{eq:convexProg:deterministic} often exhibits relatively large fluctuations. These fluctuations are unfavorable in practice since 1) servers cannot be turned on/off too frequently, and 2) the drastic change in server inlet temperature may compromise system reliability. In this section, we propose adding a regularization term to the objective function in the optimization problem~\eqref{eq:convexProg:deterministic} for suppressing undesirable fluctuations. Furthermore, we show that this modification preserves the convexity of MPC control problems. 

For a positive vector~$\xi$ of length $n_\xi$, define the posynomial
\begin{equation*}
V(\xi) = \sum_{i=1}^{n_{\xi}-1} \left(\frac{\xi_i}{\xi_{i+1}} + \frac{\xi_{i+1}}{\xi_{i}} \right).
\end{equation*}
This function measures the amount of fluctuations in the sequence~$\{\xi_1, \dotsc, \xi_{n_\xi}\}$ because the function $V$ attains its minimum if and only if $\xi_1 = \cdots = \xi_{n_\xi}$, i.e., when the sequence has no fluctuation. We then formulate an alternative finite-horizon optimal control problem as follows:
\begin{equation}\label{eq:mpc:deterministic:var}
\begin{aligned}
\minimize_{T_c,\,m_1,\,\dotsc,\,m_J,\,\gamma}\ \ \,& \begin{multlined}[t]
\gamma + w_T V(T_{c}(\tau:\tau+t_h))  \\+\sum_{j=1}^J w_j V(m_{j}(\tau:\tau+t_h))
\end{multlined}
\\
\st\ \ \ & \mbox{\eqref{eq:const:T_c}--\eqref{eq:cpuTemperatureConst}, \eqref{eq:performance_gamma},}
\end{aligned}    
\end{equation}
where $f(\tau:\tau+t_h)$ denotes the real vector $[f(\tau)\ \cdots \ f(\tau+t_h)]^\top$ for a function $f$. Adding the term $w_T V(T_{c}(\tau:\tau+t_h)) +\sum_{j=1}^J w_j V(m_{j}(\tau:\tau+t_h))$ to the original objective, we expect that the resulting control inputs show less fluctuations than the one from the previous optimization problem~\eqref{eq:mpc:deterministic}.

The next theorem shows that the optimization problem~\eqref{eq:mpc:deterministic:var} can be transformed to a convex optimization problem. The proof of the theorem is omitted due to limitations of space. 

\begin{thm}
Assume \eqref{eq:asm:T>11}. 
Then, the solution of the optimization problem~\eqref{eq:mpc:deterministic:var} is given by \eqref{eq:optimalInputRepresentation}, where 
$x^\star = \{x^\star(t)\}_{t=\tau}^{\tau+t_h}$ and $y_j^\star = \{y_j^\star(t)\}_{t=\tau}^{\tau+t_h}$ ($j=1, \dotsc, J$) solve the convex optimization problem 
\begin{align*}
\minimize_{x,\,y_1,\,\dotsc,\,y_J,\,\Gamma}\ \ \ & \begin{multlined}[t]
\log\Bigl(e^\Gamma + w_T V(e^{x(\tau:\tau+t_h)})  \\+\sum_{j=1}^J w_j V(e^{y_j(\tau:\tau+t_h)})\Bigr)
\end{multlined} %\label{eq:convexProg:deterministic:varobj:}
\\
\st\ \ \ & 
%\eqref{eq:exponentialInputSeq},\,\eqref{eq:x_bound},\,\eqref{eq:exp:delayConst},\,\eqref{eq:exp:cpuTempConst},\,\eqref{eq:lessthanGamma}. 
\mbox{\eqref{eq:exponentialInputSeq}--\eqref{eq:lessthanGamma}}. 
\end{align*}
\end{thm}

%\begin{pf}
%The convexity of the constraints \eqref{eq:exponentialInputSeq}, \eqref{eq:x_bound}, \eqref{eq:exp:cpuTempConst}, \eqref{eq:exp:delayConst}, and \eqref{eq:lessthanGamma} is already shown in the proof of Theorem~\ref{thm:deterministic}. Let us show the convexity of the objective function~\eqref{eq:convexProg:deterministic:varobj:}. Since the sum of log-convex function is convex, it is sufficient to show the convexity of the mappings
%%\begin{equation*}
%$\mathbb{R}^{t_h+1} \to \mathbb{R} \colon x(\tau:\tau+t_h) \mapsto \log V(e^{x(\tau:\tau+t_h)})$
%%\end{equation*}
%and 
%$
%%\begin{equation*}
%\mathbb{R}^{t_h+1} \to \mathbb{R} \colon y_j(\tau:\tau+t_h) \mapsto \log V(e^{y_j(\tau:\tau+t_h)})$.
%%\end{equation*}
%The convexity of these mappings follows from Lemma~\ref{lem:llconvexity} and the fact that the function $V$ is a posynomial.
%\end{pf}

\section{Scenario-based MPC}\label{sec:stochastic}

In practice, precise predictions on the user request rate (required in Assumption~\ref{asm:deterministic}) may not be available to the operator of the data center. To overcome this limitation, in this section we introduce a scenario-based MPC formalism. Specifically, we model the number of requests~$L_j(t)$ as random variables and place the following assumption in this section.

\begin{assumption}\label{asm:predictDistribution}
There exists a positive integer $t_h$ such that, at each time~$\tau$, a set of~$N$ independent sample paths
\begin{equation}\label{eq:realizations}
\{L_j^{(k)}(t)\}_{1\leq j\leq J,\,\tau\leq t\leq \tau+t_h}\quad (k=1, \dotsc, N)
\end{equation}
of the stochastic process~$L_j$ are available. 
\end{assumption}

\begin{remark}
Assumption~\ref{asm:predictDistribution} is not restrictive under our problem setting for the following reason. Typically, the operator of the data center can use a long history of past data on user requests. Furthermore, the history of data can often be regarded as a realization of a stationary stochastic process. In this case, the operator can easily obtain statistically-sound sample paths~\eqref{eq:realizations} from the past data. An implementation of this procedure is discussed in Section~\ref{sec:simulation}. 
\end{remark}

Let $E^{(k)}(\cdot; u_\tau)$ and $T_{\cpu, j}^{(k)}(\cdot; u_\tau)$ denote the trajectories of the system when the user requests rates at times~$t = \tau, \dotsc, \tau + t_h$ are those in~\eqref{eq:realizations} and the control input~$u_\tau$ of the form~\eqref{eq:def:realInputSeq} is applied. Then, we can formulate the following scenario-based optimization problem: 
\begin{align}
\minimize_{T_c,\,m_1,\,\dotsc,\,m_J,\,\gamma}\ \ \,& \begin{multlined}[t]
\gamma + w_T V(T_{c}(\tau:\tau+t_h))  \\+\sum_{j=1}^J w_j V(m_{j}(\tau:\tau+t_h))
\end{multlined}
\notag 
\\
\st\ \ \ & \eqref{eq:exponentialInputSeq},\,\eqref{eq:x_bound}, 
\notag 
\\
&
\log T_{\cpu, j}^{(k)}(t; u_\tau) \leq \log \bar T_{\cpu, j}, \label{eq:convexProg:stochastic}
\\
& y_j(t) \geq \log(\bar D_j^{-1} + L_j^{(k)}(t)) - \log \mu_j, 
\notag 
\\
& \log \sum_{t=\tau}^{\tau + t_h}  E^{(k)}(t;  u_\tau) \leq \Gamma. 
\notag 
\end{align}

In order to evaluate the reliability of the control inputs obtained by iteratively solving the optimization problem~\eqref{eq:convexProg:stochastic}, we introduce the following definition adopted from~\citet{Campi2008}. For simplicity of presentations, we hereafter focus on the case of $w_T = w_j = 0$ (i.e., no penalty on fluctuations). 

\begin{definition}[\citet{Campi2008}]
The \emph{satisfaction probability} of a control input $u_\tau \in (\mathbb{R}_+^{t_h+1})^{J+1}$ and a performance level~$\gamma$, denoted by $s(u_\tau, \gamma)$, is defined as the probability that inequalities~\eqref{eq:delayConst}, \eqref{eq:cpuTemperatureConst}, and~\eqref{eq:performance_gamma} hold true for all $\tau \leq t\leq  \tau+t_h$.
\end{definition}

The following theorem allows us to estimate the satisfaction probability. The proof of the theorem follows from the convexity of the optimization problem~\eqref{eq:convexProg:stochastic} (after the logarithmic transformations used in Section~\ref{sec:deterministic}) and Theorem~1~in \cite{Campi2008}. The details are omitted due to limitations of space.

\begin{thm}
Assume that inequality~\eqref{eq:asm:T>11} holds true. Let $x^\star = \{x^\star(t)\}_{t=\tau}^{\tau+t_h}$, $y_j^\star = \{y_j^\star(t)\}_{t=\tau}^{\tau+t_h}$, and $\Gamma^\star$ be the solutions of the convex optimization problem \eqref{eq:convexProg:stochastic}. Let $u^\star_\tau = \{e^{x^*(t)}, e^{y_1^*(t)}, \dotsc, e^{y_J^*(t)} \colon 1\leq j\leq J,\,\tau\leq t\leq \tau+t_h\}$. Then, $\Pr(s(u_\tau^\star, \log \Gamma^\star) > \epsilon) \leq \sum_{i=0}^{(t_h+1)J}\binom{N}{i} \epsilon^i (1-\epsilon)^{N-i}$.
\end{thm}

\newcommand{\figureDate}{D05-Apr-2018_19-52-04_S05-Apr-2018_19-48-48}
\newcommand{\folderName}{../../simulations/figures/\figureDate}
\renewcommand{\folderName}{}

\section{Numerical Simulations}\label{sec:simulation}

In this section, we illustrate the effectiveness of the theoretical results obtained in the previous sections by using the real-world Google Workload Data Trace~\citep{Reiss2011} of~$J=3$ users. We pre-process the data so that the average workload arrival rate is calculated per minute. The first $\num[group-separator={,}]{30000}$ arrival rate samples are used to train our scenario-based MPC model. To evaluate the performance of the proposed MPC policy (Section~\ref{sec:deterministic}) and scenario-based MPC policy (Section~\ref{sec:stochastic}), we introduce an Offline Optimal Static (OOS) policy, which uses a constant CRAC supplying air temperature in all decision epochs to minimize the overall energy consumption while enforcing constraints \eqref{eq:const:T_c}, \eqref{eq:delayConst}, and \eqref{eq:cpuTemperatureConst}, provided that the controller has full system information in the future.

Note that the constant temperature policy is widely used in the industry, and OOS presumably has better performance than the conventional constant temperature policy since OOS is equipped with future system information. The deterministic MPC policy requires, at each time $t= \tau$, correct predictions of $L_j(t)$ from $t=\tau$ to $t = \tau + t_h$.  On the other hand, at each time, the scenario-based MPC policy only needs to know the request rate at the current time as well as the past history of the user request rates.

To implement our scenario-based MPC policy, we find $N=100$ different sequences $\{L_{j}(t)\}_{t = t_{0k}}^{t_{0k} + t_h}$ ($k=1, \dotsc, N$) for each time $t=\tau$ and user $j$ such that $L_j(t_{0k}) = L_j(\tau)$ from the past data, and use the sequences as the sample paths~\eqref{eq:realizations} required by Assumption~\ref{asm:predictDistribution}. If there are not enough number of sequences satisfying $L_j(t_{0k}) = L_j(\tau)$, we instead use a sequence such that $L_j(t_{0k})$ is as close to $L_j(\tau)$ as possible.

\subsection{Result Analysis}

We set $t_h = 5$, $a_1 = 10$, $a_2 = 1$, $\alpha_j = 0.05$, $\beta_j = 0.95$, $\varsigma = 1.5$, and  $\mu_j = 1$. We assume that the cooling temperature can be tuned within the interval $18 \leq T_c(t)\leq 27$, which satisfies the assumption~\eqref{eq:asm:T>11}. The maximum CPU temperature and the upper bound of response times are set as $\bar T_{\cpu, j} = 80$ and~$\bar D_j = 0.05$, respectively. We use the initial value~$T_{\cpu, j}(t_0) = 27$. We let $t_0 = \num[group-separator={,}]{30000}$ and~$t_f = \num[group-separator={,}]{30200}$. For simplicity, the weight, $w_m$, for penalizing the fluctuations in the number of servers is set to be zero. 

\begin{figure}[tb]
\centering
\includegraphics[width=.95\linewidth,trim={0mm 0mm 0mm 1mm},clip]{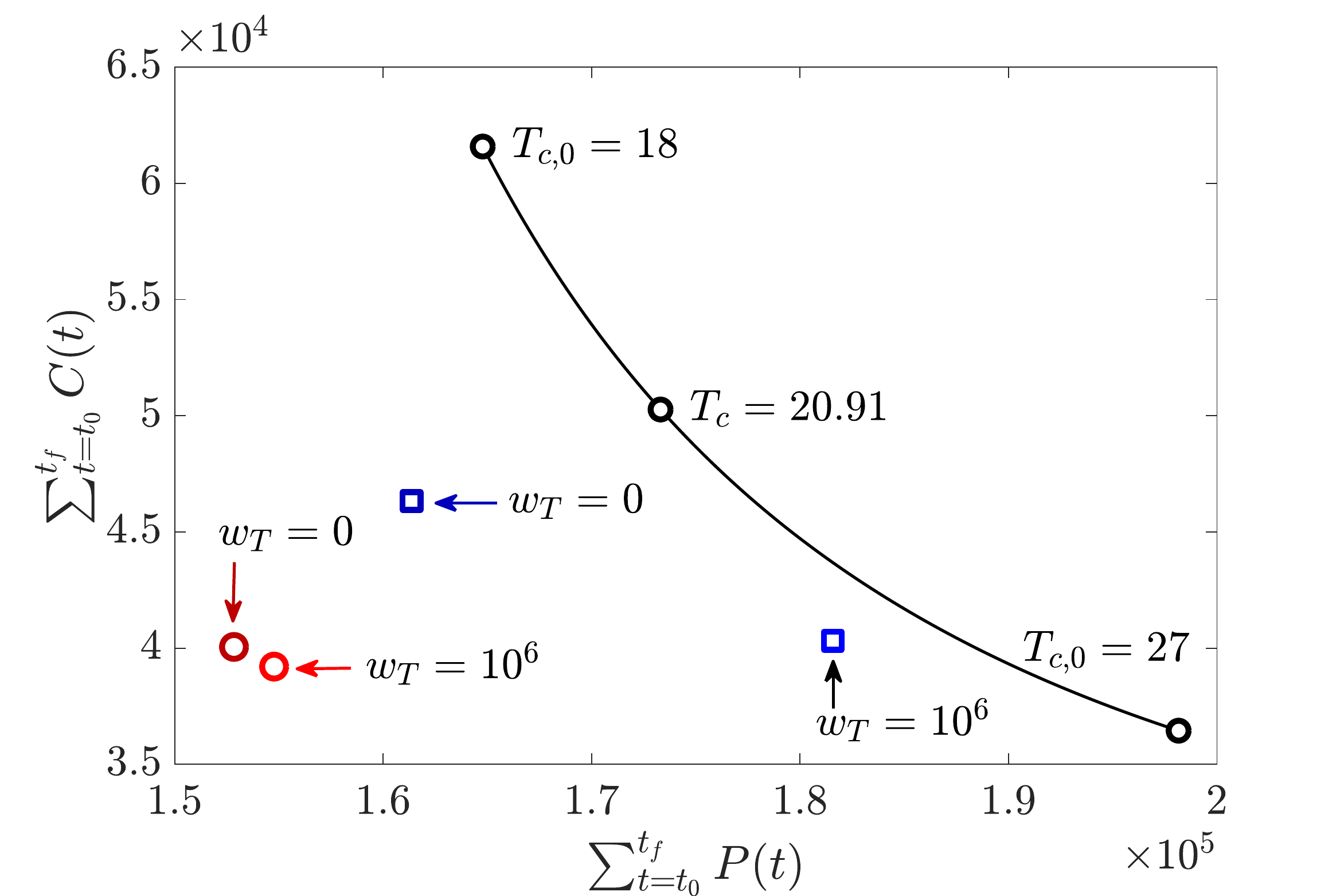}
\caption{Total IT and cooling energy. Circles: MPC policies, Squares: Scenario-based MPC policies, Solid line: Offline Optimal Static (OOS) policy.}
\label{fig:deterministic:comparisonEnergy}
\end{figure}

\begin{figure}[tb]
\centering
\includegraphics[width=.95\linewidth,trim={0mm 3mm 0mm 5mm},clip]{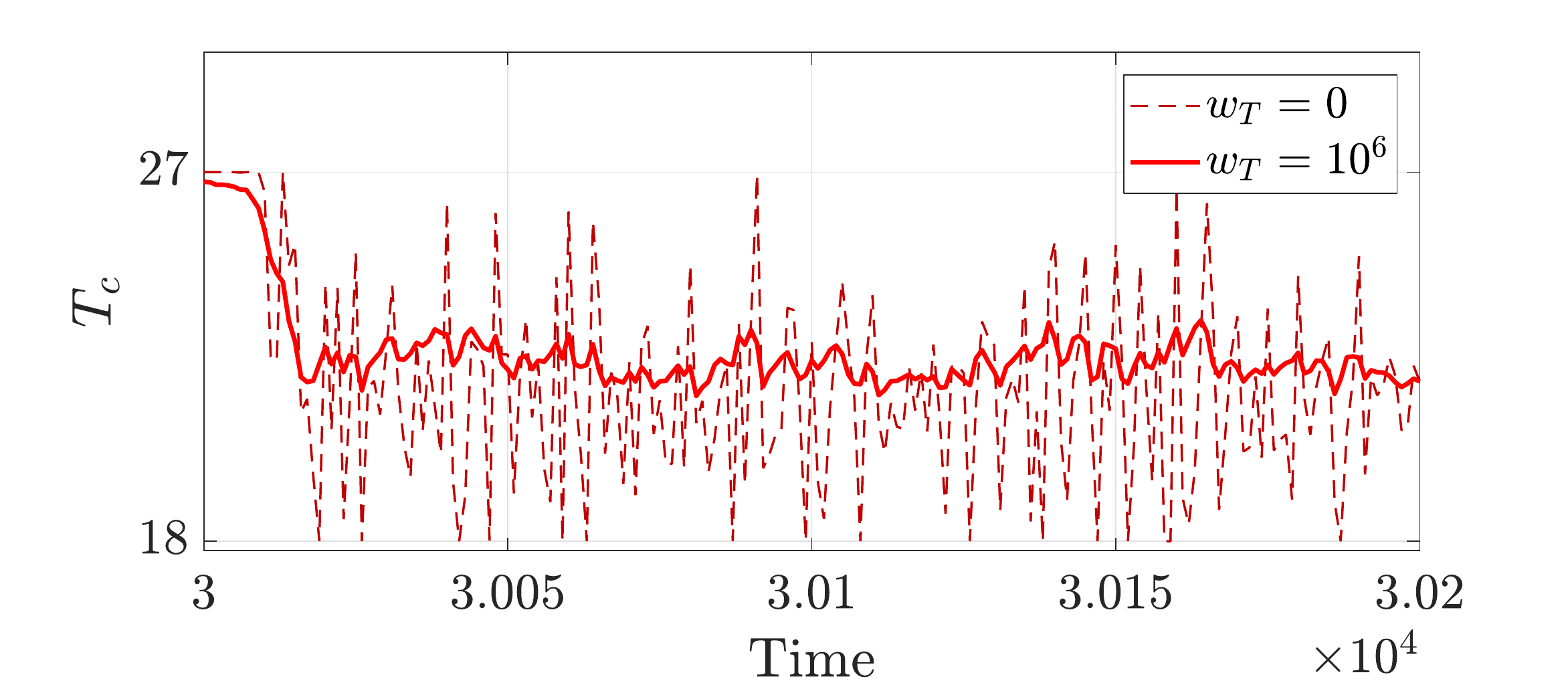}
\caption{Cooling temperatures with MPC policies.}
\label{fig:deterministic:temperatures}
\vspace{3mm}
\centering
\includegraphics[width=.95\linewidth,trim={0mm 3mm 0mm 5mm},clip]{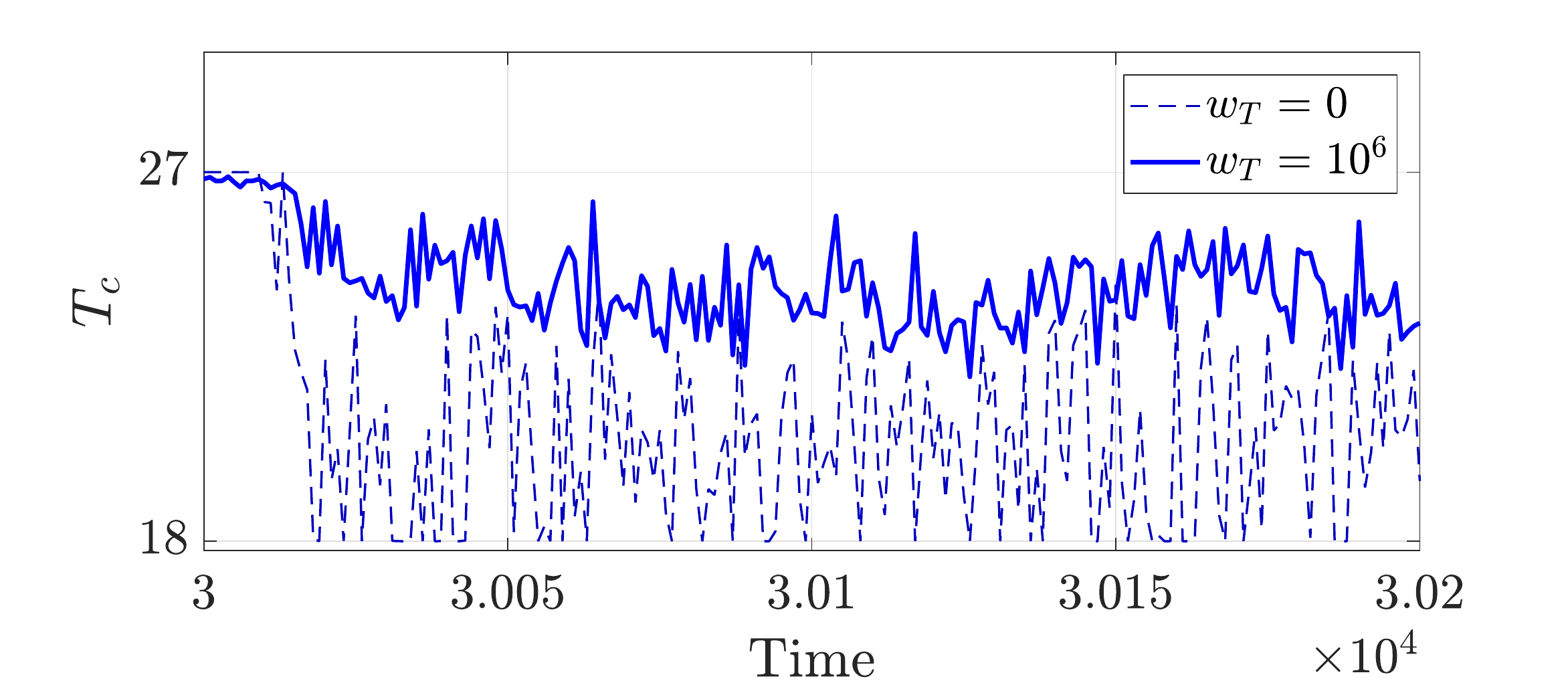}
\caption{Cooling temperatures with scenario-based MPC policies.}
\label{fig:stochastic:temperatures}
\end{figure}

In Fig.~\ref{fig:deterministic:comparisonEnergy}, we compare the energy consumption for the three strategies. We see that the MPC policy without penalization ($w_T = 0$) achieves the best performance. The introduction of the penalization term in~\eqref{eq:mpc:deterministic:var} with $w_T= 10^6$ only slightly deteriorates the performance by $0.56\%$. The scenario-based MPC without penalization consumes $7.67\%$ more energy compared to its MPC counterpart, and the penalization further mildly degrades the performance by $7\%$. However, the scenario-based approach still saves $0.8\%$ overall energy compared to OOS, which achieves the lowest energy consumption when $T_c=20.9^\circ$C. Note that OOS requires full future information and is the opportunistic lower bound of practical static policies. We here also emphasize that the scenario-based policy is the most realistic among other policies because it does not require predictions, which is hard to be precise in practice.

Figs.~\ref{fig:deterministic:temperatures} and \ref{fig:stochastic:temperatures} show the instant CRAC supplying air temperature $T_c$ for MPC and scenario-based MPC, respectively. If there is no penalization term, $T_c$ in MPC is generally $1^\circ$C higher than scenario-based MPC. However, if the penalization term is presented, $T_c$ in MPC is $1.7^\circ$C lower than scenario-based MPC. Since a higher $T_c$ increases the CRAC efficiency, we can conclude that scenario-based MPC with penalization has the potential to save the cooling energy. In addition, it is shown that the penalization greatly suppresses fluctuations in $T_c$. More specifically, the variance of $T_c$ drops from $5.8$ to $1.1$ for MPC, and from $6.0$ to $1.2$ for scenario-based MPC.

\begin{figure}[tb]
\centering
\includegraphics[width=.95\linewidth,trim={-1mm 15mm 0mm 18.5mm},clip]{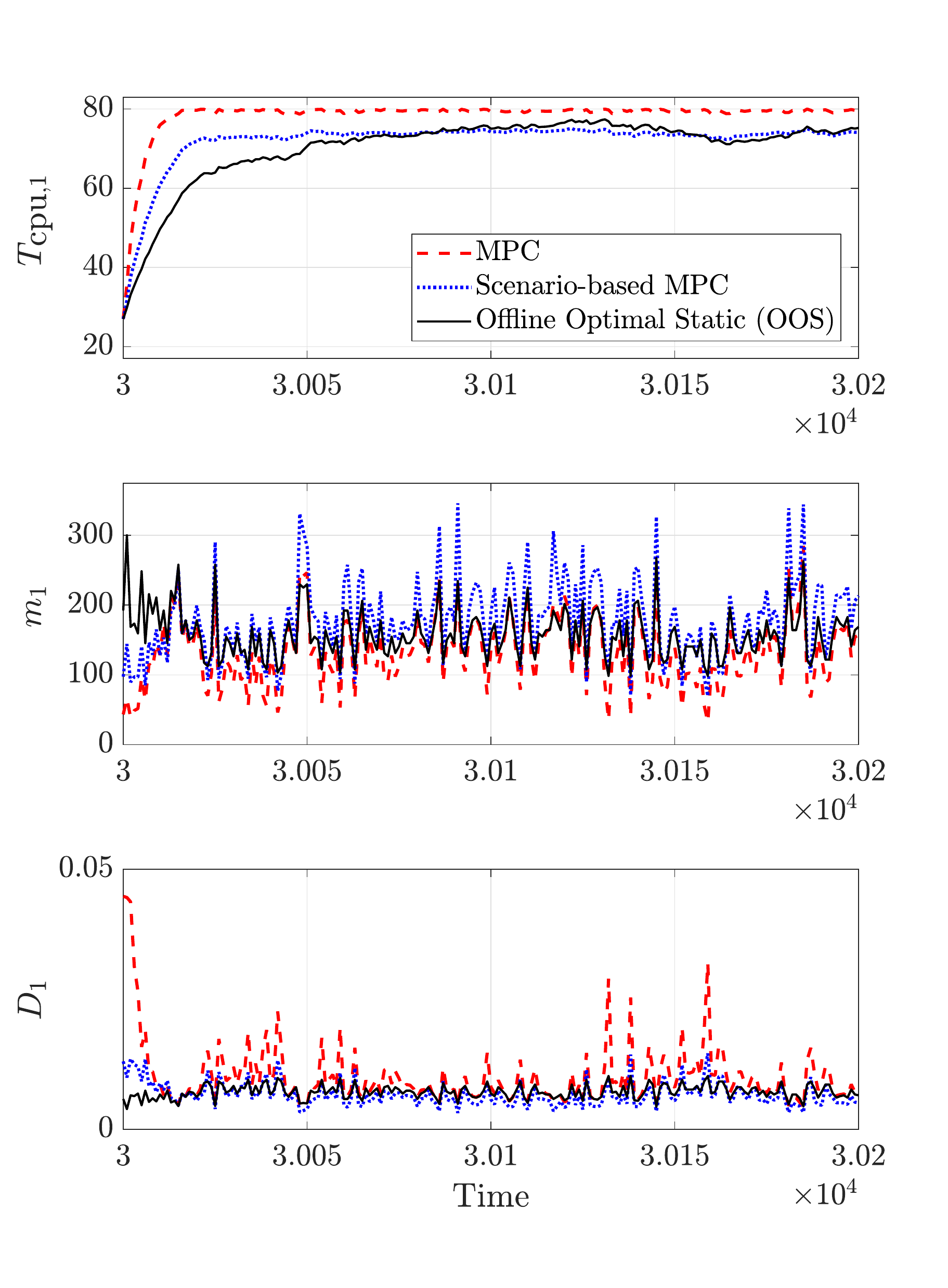}
\caption{Servers allocated to user 1. From top to bottom: the CPU temperature, the number of the servers, and the response time.}
\label{fig:deterministic:sever1}
\end{figure}

Finally, we illustrate the trajectories of the CPU temperature, the number of  servers, and the response time for the first cluster of servers in Fig.~\ref{fig:deterministic:sever1}. We can see that OOS generally achieves the lowest CPU temperature, since the constant CRAC supplying air setpoint has to be low enough to guarantee the thermal constraint \eqref{eq:cpuTemperatureConst} in the worst case. In contrast, the CPU temperature in deterministic MPC is very close to the threshold $\bar{T}_{\cpu,j}$ since the supplying air temperature can be dynamically tuned. Scenario-based MPC uses more servers than basic MPC due to the lack of future system information. The delay constraint \eqref{eq:delayConst} is satisfied for all three policies.

\section{Conclusions}

In this paper, we have studied the problem of minimizing the energy consumption in data centers with cold aisle containment. We have proposed standard and scenario-based MPC policies to solve the problem. Since our formulation allows the sequences of optimal control problems to be reduced to convex optimization problems, we can efficiently and uniquely identify control inputs with a guaranteed optimality. With several numerical simulations, we have confirmed the effectiveness of the proposed policies.


%\section*{ACKNOWLEDGMENT}

%%

%\bibliography{Mendeley.bib}

\end{document}